\documentclass[12pt]{article}

\usepackage{pstricks}
\usepackage[dvips]{graphicx} 

\usepackage{amsmath}
\usepackage{amssymb}

      \newcommand {\al}   {\alpha}          \newcommand {\bt}  {\beta}
                
      \newcommand {\del}  {\delta}          
              \newcommand {\ve}   {\varepsilon}
                 \newcommand {\vphi} {\varphi}
      \newcommand {\lam}  {\lambda}         
      
      \newcommand {\om}   {\omega}          \newcommand {\Om}  {\Omega}
      \newcommand {\pl}   {\partial}

      \newcommand {\RRR}  {{\mathbb R}}     \newcommand {\RRRR}  {{\cal R}}
              \newcommand {\LLL}  {{\cal L}}
              
           \newcommand {\NNN}  {{\cal N}}
      \newcommand {\III}  {{\cal I}}        
      \newcommand {\FFF}  {{\cal F}}        
              
      \newcommand {\interval}  {[-\pi/2,\, \pi/2]}

\author{Alexander Plakhov\thanks{Department of Mathematics, Aveiro
University, Aveiro 3810, Portugal} \and Paulo
Gouveia\thanks{Techn. and Manag. School, Bragan\c{c}a Polytechnic Institute,
5301 Bragan\c{c}a, Portugal}}

\title{Problems of maximal mean resistance\\ on the plane}

\date{}

\begin{document}

\maketitle

\begin{abstract}
A two-dimensional body moves through a rarefied medium; the
collisions of the medium particles with the body are absolutely
elastic. The body performs both translational and slow rotational
motion. It is required to select the body, from a given class of
bodies, such that the average force of resistance of the medium to
its motion is maximal.

There are presented numerical and analytical results concerning this
problem. In particular, the maximum resistance in the class of
bodies contained in a convex body $K$ is proved to be 1.5 times
resistance of $K$. The maximum is attained on a sequence of bodies
with very complicated boundary. The numerical study was made for
somewhat more restricted classes of bodies. The obtained values of
resistance are slightly lower, but the boundary of obtained bodies
is much simpler, as compared to the analytical solutions.
\end{abstract}

\begin{quote}
\end{quote}

\begin{quote}
\end{quote}

\begin{quote}
{\small {\bf Running title:} Bodies of maximal resistance}
\end{quote}

\section{Introduction}

Consider a homogeneous medium of point particles at rest in
Euclidean space $\RRR^d$, and a body moving forward through this
medium. The medium is highly rarefied, so that mutual interaction of
the particles is neglected. The interaction of the particles with
the body is absolutely elastic. It is required to find a shape of
the body that minimizes or maximizes resistance of the medium to its
motion.

When thinking of this kind of problems, one can have in mind an
artificial satellite of the Earth moving on a relatively low (say,
$100 \div 200$ km) orbit; one has to {\it minimize} resistance of
the rest of atmosphere, or a solar sail; then one has to {\it
maximize} the pressure of the flux of solar photons on the sail.

In order to specify the problem, one has to describe the body's
motion as well as to define the class of admissible bodies. In
classes of {\it convex} bodies with {\it translational} motion, the
minimization problem has been extensively studied. Newton \cite{N}
obtained the solution in the class of (three-dimensional) convex
{\it axially symmetric} bodies of fixed length and width. Since
1993, there have been obtained many interesting results in classes
of convex {\it non-symmetric} bodies \cite{BK}-\cite{LO}.

Note that resistance can be written in the form $R[f] =
\int\!\!\int_D (1 + |\nabla f|^2)^{-1}\, dx\, dy$, where the
function $z = f(x,y)$,\, $(x,y) \in D$ describes the front part of a
convex body. Here the (orthogonal) coordinates are chosen in such a
way that the body's velocity equals $(0, 0, 1)$. Thus, the problem
amounts to minimization of the functional $R[f]$.

This approach is not valid, as applied to {\it nonconvex} bodies. (A
body is a bounded connected set with piecewise smooth boundary.) The
reason is that the above formula is not true if particles can hit
the body more than once. In fact, there is no simple analytic
formula for resistance in the nonconvex case.

In general it is not easy to calculate resistance even for
(nonconvex) bodies with well-behaved boundary; however, usually one
can construct a minimizing body or a minimizing sequence of bodies.
In the three-dimensional case infimum of resistance is typically
equal to zero \cite{P1,P2}\footnote{This result is obtained for
classes of bodies of fixed length and width \cite{P2} and for
classes of bodies containing a bounded set and being contained in
its $\ve$-neighborhood \cite{P1}.}. On the contrary, in the
two-dimensional case infimum of resistance is positive and usually
can be explicitly found \cite{P2}. In higher dimensions, $d
> 3$, the answer is the same as in the three-dimensional case:
infimum of resistance equals zero. Note in passing that the problem
of {\it maximal} resistance admits a very simple solution: the front
part of the body's surface should be orthogonal to the direction of
motion, or should be composed of pieces orthogonal to this
direction.

It is also interesting to consider rotational motion of the body.
Imagine an artificial satellite without orientation control; one can
expect that in the course of motion, it will perform a (perhaps very
slow) rotation. The problem of {\it minimal mean resistance}  for
non-convex {\it rotating} bodies in two dimensions was considered in
\cite{average-04}; it was proved that the gain in resistance, as
compared with the convex case, is smaller than $1.22\%$.

Here we study the problem of {\it maximal mean resistance} for {\it
rotating} bodies in two dimensions, $d = 2$. This problem is far
from being trivial, contrary to the case of purely translational
motion. To see it, consider a unit disk on the plane, which is
moving forward and at the same time slowly (and uniformly) rotating.
Denote the disk by $K_1$. ``Cut off''{} a small portion of the disk
contained in the $\ve$-neighborhood of $\pl K_1$ ($\ve \ll 1$); the
resulting set $B$ is such that $B \subset K_1 \subset \NNN_\ve (B)$
(here $\NNN_\ve$ designates $\ve$-neighborhood). The question is:
how large can the increase of resistance
$\,\frac{\text{Resistance}(B)}{\text{Resistance}(K_1)}\,$ be? Some
estimates can be made immediately. Firstly, it cannot exceed
$1{.}5$. This hypothetical maximal increase is achieved if the
velocity of a reflected particle is always opposite to the incidence
velocity, $v^+ = -v$; in this case the momentum transmitted by the
particle to the body is maximal. Next, if the circumference $\pl
K_1$ is partitioned into several small arcs and each arc is
substituted with a pair of legs of a right isosceles triangle
contained in $K_1$ (the resulting body is shown on Fig.\,(a) below)
then resistance increases approximately $\sqrt 2$ times. More
precisely, if the length of each arc is $2\ve$ then the resistance
increase is $\sqrt 2\, {\sin\ve}/{\ve} \approx \sqrt 2 (1 -
\ve^2/6)$; for the proof see Appendix 1.

Another example is the body obtained by making deep and narrow
rectangular ``hollows''{} on the boundary of $K_1$; see Fig.\,(b).
Put the depth of any hollows to be $\ve$, width $\ve^2$, and the
distance between neighboring hollows, $\ve^3$. Then approximately
one half of the particles incident on a hollow get out with the
velocity opposite to the initial one, and the rest of the particles
get out with the velocity symmetric to the initial one with respect
to the smaller side of the rectangle. Resistance of the body is,
approximately, the arithmetic mean of the disk resistance and the
(hypothetical) maximal resistance (1.5 times the disk resistance);
that is, the resistance increase is $1.25 + o(1)$,\, $\ve \to 0$.
\begin{center}
\begin{tabular}{c c}
\includegraphics[width=0.35\columnwidth]{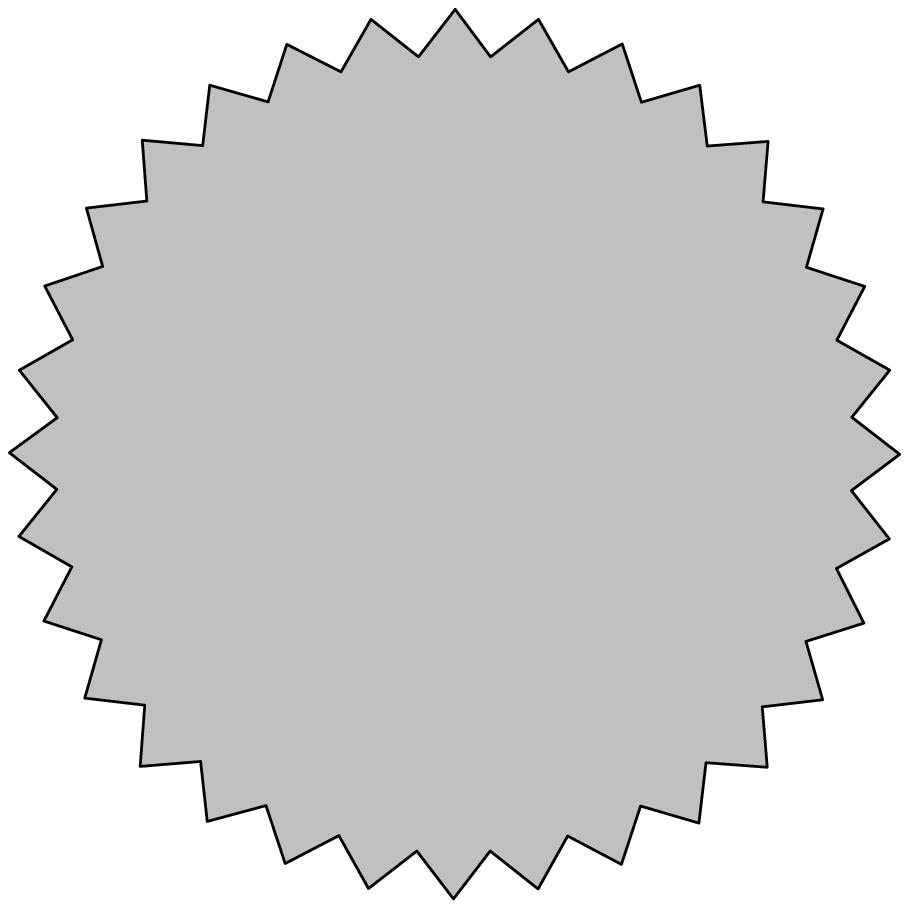}
&
\includegraphics[width=0.35\columnwidth]{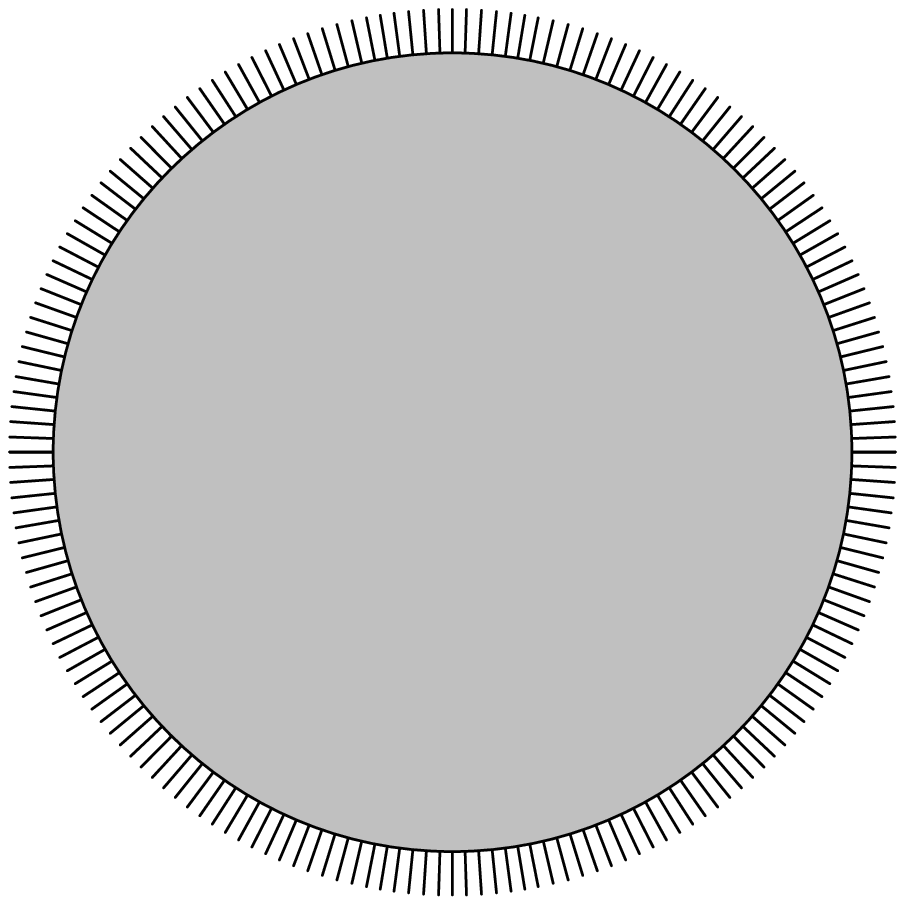} \\
(a)&(b)\\
\end{tabular}
\end{center}

In this paper, the resistance maximization problem is studied

(\~A) in the class of planar sets of the form $r \le 1 - \ve
f(\theta/\ve)$,\, $0 < \ve \ll 1$,\, $f \ge 0$ in polar coordinates
$r$,\, $\theta$ and

(\~B) in the class of sets contained in a fixed two-dimensional
convex body.

These problems still have to be rigorously stated; it is made in
Section 2, and the corresponding reformulated problems take the
names (A) and (B). The restricted problem (A) seemed to be more
amenable to numerical study and was examined first. It was not
completely solved; nevertheless, we present here some numerical
results. We believe that they are of interest, since they allow one
to reach values up to 1.446 (which is rather close to the upper
bound of resistance) by using relatively simple geometric shapes.

When working on the restricted problem, there was gained experience
that eventually allowed one to solve the problem (B). The answer
here is 1.5; the corresponding maximizing sets have much more
complicated boundary than those used in (A).

The paper is organized as follows. The mathematical formulation of
the problems is given in Section 2. The problem (A) is studied
numerically in Section 3, and the problem (B) is solved analytically
in Section 4. The obtained results are discussed in Conclusions.
Finally, some auxiliary formulas related to resistance of zigzag
shapes are derived in Appendices 1 and 2.

\section{Statement of the problem}

Let $B \subset \RRR^2$ be a bounded connected set  with piecewise
smooth boundary; consider the billiard in $\RRR^2 \setminus B$.
Consider a billiard particle that initially moves freely, then makes
several (at least one) reflections at regular points of $\pl B$, and
finally, moves freely again. Denote by $\text{conv}B$ the convex
hull of $B$.

The trajectory of the particle intersects $\pl(\text{conv}B)$ twice:
when getting in the set $\text{conv}B$ and when getting out of it.
Note that if the point of first intersection belongs to $I_0 :=
\pl(\text{conv}B) \cap \pl B$, then the two points of intersection
coincide. Introduce the natural parametrization of
$\pl(\text{conv}B)$ by the parameter $\xi \in [0,\, L]$, where $L =
|\pl(\text{conv}B)|$ is length of the curve $\pl(\text{conv}B)$.
 \vspace{0mm}

Let $\xi$ and $\xi^+$ be the first and second intersection points,
and let $v$ and $v^+$ be the particle velocity at these points,
respectively. Denote by $\langle \cdot\,, \cdot \rangle$ the scalar
product. Let $n_\xi$ be the outward unit normal vector to
$\pl(\text{conv}B)$ at the point corresponding to $\xi$. For a
vector $w$ such that $\langle w, n_\xi \rangle \ge 0$, let us agree
that the angle between $n_\xi$ and $w$ is counted from $n_\xi$ to
$w$ clockwise or counterclockwise; in the first case it is positive,
and in the second case, negative. Thus, the angle varies in the
interval $\interval$. Denote by $\vphi$ and $\vphi^+$ the angles the
vectors $n_\xi$ and $n_{\xi^+}$ form with $-v$ and $v^+$,
respectively. (Note that one always has $n_\xi = n_{\xi^+}$.)
Thereby, the one-to-one mapping $T_B: (\vphi, \xi) \to (\vphi^+,
\xi^+)$ is determined. It is defined and takes values on a full
measure subset of $\interval \times [0,\, L]$. Moreover, the
following holds true.

(i) $T_B$ preserves the measure $\mu$ given by $d\mu(\vphi,\xi) =
\cos\vphi\, d\vphi d\xi$;

(ii) $T_B^{-1} = T_B$.\\
These relations follow from the measure preserving property and from
time-reversibility of billiard dynamics; see \cite{average-04} for
more details.
 \vspace{0mm}

\vspace*{70mm}

\scalebox{0.7}{
 \pscurve[fillstyle=solid,
fillcolor=lightgray](3,8)(2.8,7.7)(2.5,6.5)(1.7,5.4)(0.5,4.7)(0.1,4)(0,2.5)(1,1)(3,0.2)(4.5,2)(5,1.4)(5.5,1.4)(5.5,0.8)(6,0.2)(9,1.25)(11.4,3)
 (12.3,4.5)(12.4,5.35)(12,5.5)(11,5.2)(9.5,5.2)(8.3,6)(7.9,7)(8,8)
 \pscurve[fillstyle=solid,
fillcolor=white](3,8)(3.2,8)(4,6.9)(5,6.4)(5.6,7)(5.9,7.2)(6.5,6.3)(7.6,7)(8,8)
  \psline[linestyle=dashed,linewidth=0.4pt](2.7,0.16)(6.4,0.14)
 \psline[linestyle=dashed,linewidth=0.4pt](3,8)(8,8)
 \psline[linestyle=dashed,linewidth=0.4pt](2.98,8)(0.35,4.55)
 \psline[linestyle=dashed,linewidth=0.4pt](12.3,5.45)(8,8)
 \psline[linewidth=1.2pt,arrows=->,arrowscale=1.5](8.3,9.5)(8.3,6)(11,5.2)(12.6,6.4)
 \psline[linewidth=1.2pt,arrows=->,arrowscale=1.5](8.3,9.55)(8.3,8.5)
 \psecurve[fillstyle=solid,
fillcolor=white](7,3)(5.5,3.5)(7,4.4)(7.5,4.6)(8,3.5)(8.1,3.3)(7.7,3)(7,2.5)(6.3,3)(5.5,3.5)(7,4.4)
 \rput(8.55,9){\Large $v$}
 \rput(13,6.3){\Large $v^+$}
 \rput(8.55,8){\Large $x$}
 \rput(11.8,6.2){\Large $x^+$}
 \psdots[dotsize=3pt](8.3,7.82)(11.76,5.77)
    \rput(9.25,0.9){\Large $I_0$}
    \rput(10.4,7.1){\Large $I_3$}
    \rput(4,-0.2){\large $I_4$}
    \rput(9.8,6.2){$\Om_3$}
   \rput(4.7,0.8){\Large $\Om_4$}
    \rput(7,3.6){\Large $\Om_5$}
  \rput(5.8,8.35){\Large $I_2$}
     \rput(4.8,7.2){\large $\Om_2$}
    \rput(1.7,6.9){\Large $I_1$}
     \rput(1.85,5.95){$\Om_1$}
  }
   \rput(2.5,2.5){\Huge $B$}

 \vspace{10mm}

Suppose now that the center of mass of the body moves forward at the
velocity $e_2 = (0, 1)$ and the body rotates around the center of
mass with a small angular velocity $\om \ll 1$. Thus, each
individual particle interacts with the body in the same way as if it
there were no rotation. Resistance of the medium is a periodic
vector-valued function of time, $R_B(t)$, with the period $T =
2\pi/\om$. In order to derive the formula for the mean value of
resistance $\RRRR(B) = \frac 1T\, \int_0^T R_B(t)\, dt$, consider a
reference system moving forward at the velocity $e_2$. In this
reference system, the body rotates around a fixed point, and there
is a flux of particles of velocity $-e_2$ incident on the body. Each
particle transmits to the body a momentum proportional to $v - v^+$,
where $v = -e_2 = (0, -1)$ is the initial velocity of the particle,
$v^+ = (\sin(\vphi^+ - \vphi), \cos(\vphi^+ - \vphi))$ is its final
velocity, $\vphi = \vphi_{0i} + \om t$,\, $t$ is the moment of the
first intersection of the particle with $\pl(\text{conv}B)$,\,
$\vphi^+ = \vphi^+_{(B)}(\vphi,\xi)$; here $\xi$ means the point of
the first intersection of the particle with $\pl(\text{conv}B)$ and
$\vphi_{(B)}^+(\vphi,\xi)$ is the first component of the mapping
$T_B$. The mean resistance $\RRRR(B)$ is the sum of all momenta
transmitted to the body in a time interval of length $T$, divided by
$T$, that is,
\begin{equation}\label{1}
\RRRR(B) = -c \int_{-\pi/2}^{\pi/2} \int_0^L \left(
\sin(\vphi_{(B)}^+(\vphi,\xi) - \vphi),\ 1 +
\cos(\vphi_{(B)}^+(\vphi,\xi) - \vphi) \right) d\mu(\vphi, \xi);
\end{equation}
the ratio $c$ is proportional to the medium density. Later on we
shall specify the appropriate value of $c$ simplifying the
subsequent formulas.

Further, changing the variables $(\vphi, \xi) \mapsto (\tilde\vphi,
\tilde\xi) = (\vphi_{(B)}^+(\vphi,\xi),\, \xi_{(B)}^+(\vphi,\xi))$
in the integral (\ref{1}) and taking into account that the Jacobian
related to this change is 1 (property (i)) and that
$\vphi_{(B)}^+(\tilde\vphi, \tilde\xi) = \vphi$ (property (ii)), one
gets
\begin{equation}\label{1'}
\RRRR(B) = -c \int_{-\pi/2}^{\pi/2} \int_0^L \left( \sin(\tilde\vphi
- \vphi_{(B)}^+(\tilde\vphi, \tilde\xi)),\ 1 + \cos(\tilde\vphi -
\vphi_{(B)}^+(\tilde\vphi, \tilde\xi)) \right) d\mu(\tilde\vphi,
\tilde\xi).
\end{equation}
Comparing (\ref{1}) and (\ref{1'}) and taking into account that sine
is odd, one concludes that the first component of $\RRRR(B)$ is
zero, that is,
\begin{equation}\label{1''}
\RRRR(B) = -c \int_{-\pi/2}^{\pi/2} \int_0^L \left( 1 +
\cos(\vphi_{(B)}^+(\vphi,\xi) - \vphi) \right) d\mu(\vphi, \xi)
\cdot e_2.
\end{equation}

Let us now reduce the formula (\ref{1''}) to the form more
convenient for computation. The curve $\pl(\text{conv}B)$ is the
union of a finite or countable family of sets $I_0$,\, $I_1$,\,
$I_2, \ldots$,
$$
\pl(\text{conv}B) = \cup_i I_i.
$$
Here $I_0 = \pl(\text{conv}B) \cap \pl B$ is the ``convex part''{} of
the boundary $\pl B$, and $\pl(\text{conv}B) \setminus \pl B$ is the
union of open intervals $I_1$,\, $I_2, \ldots$. Respectively,
$\text{conv}B$ is the union of a finite or countable family of sets
$\Om_0$,\, $\Om_1$,\, $\Om_2, \ldots$,
$$
\text{conv}B = \cup_i \Om_i,
$$
where $\Om_0 = B$ and the sets $\Om_1$,\, $\Om_2, \ldots$ are
connected components of $\text{conv}B \setminus B$ (``cavities''{} on
$B$). The enumeration is chosen in such a way that $I_i \subset
\pl\Om_i$ (see the figure above). Note that for some sets $\Om_i$
there may happen $\pl\Om_i \subset \pl B$. These sets (``interior
cavities''{} of $B$) have no influence on resistance and will be
ignored in the sequel.

Each interval $I_i\, (i \ne 0)$ corresponds to an interval $\LLL_i$
(modulo L) on the parameter set $[0,\, L]$, and $I_0$ corresponds to
$[0,\, L] \setminus \cup_{i \ne 0} \LLL_i =: \LLL_0$. Let $l_i$ be
the length of $L_i$. The sets $\interval \times \LLL_i$ are
invariant with respect to $T_B$; denote by $T_i$ the restrictions of
$T_B$ on these sets, and by $\vphi_i(\vphi, \xi)$,\, $\xi_i(\vphi,
\xi)$, the components of $T_i$. One easily sees that $T_0$ is given
by $T_0(\vphi, \xi) = (-\vphi, \xi)$. The mean resistance $\RRRR(B)$
(\ref{1''}) is the sum of partial resistances $\RRRR_i$ related to
the $i$th cavity:\, $\RRRR(B) = \sum_i \RRRR_i$, where
\begin{equation}\label{3}
\RRRR_i = -c \int_{-\pi/2}^{\pi/2} \int_{\LLL_i} \left( 1 +
\cos(\vphi_i^+(\vphi,\xi) - \vphi) \right) d\mu(\vphi, \xi) \cdot
e_2.
\end{equation}
In particular,
$$
\RRRR_0 = -c \int_{-\pi/2}^{\pi/2} \int_{\LLL_0} \left( 1 + \cos
2\vphi \right) \cos\vphi\, d\vphi\, d\xi \cdot e_2 = -\frac 83\, c\,
l_0\, e_2.
$$
Put $c = 3/8$, then one gets a slightly simplified relation:~
$\RRRR_0 = -l_0\, e_2$.

We say that a bounded set $\Om \subset \RRR^2$ with piecewise smooth
boundary is a {\it standard cavity} if it contains the interval
$\III := [0,\, 1] \times \{ 0 \}$ and is contained in the upper
half-plane $\{ (x_1, x_2):\, x_2 \ge 0 \}$, that is,
$$
\III \subset \Om \subset \{ (x_1, x_2):\, x_2 \ge 0 \}.
$$
Consider billiard in $\Om$; suppose that a billiard particle starts
from a point of $\III$, and after several reflections from $\pl\Om
\setminus \III$ returns to $\III$. Let $(\xi, 0)$ and $(\sin\vphi,
\cos\vphi)$ be the initial location and initial velocity of the
particle, and denote by $(\xi_\Om(\vphi,\xi), 0)$ and
$-(\sin\vphi_\Om(\vphi,\xi), \cos\vphi_\Om(\vphi,\xi))$, the final
point and final velocity. The so defined map $(\vphi, \xi) \mapsto
(\vphi_\Om(\vphi,\xi),\, \xi_\Om(\vphi,\xi))$ preserves the measure
$\mu$ and is defined and takes values on a full measure subset of
$\interval \times [0,\, 1]$.

   \vspace{55mm}

\scalebox{0.7}{
 \rput(5,0){
 \psline[linewidth=0.5pt,linestyle=dashed](0,0)(8,0)
 \pscurve[linewidth=0.5pt](0,0)(-1.7,3)(-1,6)(2,6.8)(4,5.8)(6,5.4)(9,5.3)(9.3,4)(8.5,3)(7.2,1.5)(8,0)
 \psline[linewidth=1pt,arrows=->,arrowscale=2](1,0)(3.15,2.75)
 \psline[linewidth=1pt](1,0)(3.15,2.75)(5.3,5.5)
 \psline[linewidth=1pt,arrows=->,arrowscale=2](5.3,5.5)(6.25,3.5)
 \psline[linewidth=1pt](6.25,3.5)(7.2,1.5)
 \psline[linewidth=1pt,arrows=->,arrowscale=2](7.2,1.5)(6.6,0.75)
 \psline[linewidth=1pt](6.6,0.75)(6,0)
 \psline[linewidth=0.5pt,linestyle=dashed](1,0)(1,1.2)
 \psline[linewidth=0.5pt,linestyle=dashed](6,0)(6,1)
  \psarc(1,0){0.8}{52}{90}
  \psarc(6,0){0.8}{52}{90}
  \rput(1,-0.3){\large $\xi$}
  \rput(6,-0.4){\large $\xi_\Om(\vphi,\xi)$}
  \rput(1.4,1.2){\large $\vphi$}
  \rput(5.9,1.4){$\vphi_\Om(\vphi,\xi)$}
}}
 \rput(1.3,2.5){\Huge $\Om$}
 \rput(5.6,-0.8){\Large $\III$}

   \vspace{15mm}

Designate
\begin{equation}\label{5}
\FFF(\Om) = \frac 38 \int_{-\pi/2}^{\pi/2} \int_0^1 \left( 1 +
\cos(\vphi_{\Om}(\vphi,\xi) - \vphi) \right) \cos\vphi\, d\xi
d\vphi;
\end{equation}
The integrand in (\ref{5}) does not exceed 2, therefore $\FFF(\Om)
\le (3/8) \int_{-\pi/2}^{\pi/2} \int_0^1 2\, \cos\vphi\, d\xi d\vphi
= 1.5$. On the other hand, denoting $\Om_\ve := \pi_\ve \Om$, where
$\pi_\ve:$ $(x_1, x_2) \mapsto (x_1, \ve x_2)$, one has that
$\vphi_{\Om_\ve}(\vphi,\xi) \to -\vphi$ as $\ve \to 0^+$, hence
$\lim_{\ve \to 0^+} \FFF(\Om_\ve) = (3/8) \int_{-\pi/2}^{\pi/2}
\int_0^1 (1 + \cos 2\vphi)\, \cos\vphi\, d\xi d\vphi = 1$. It
follows that $1 \le \sup_\Om \FFF(\Om) \le 1.5$.

Each pair $(\Om_i, I_i)$,\, $i = 1,\, 2,\ldots$ can be reduced, by a
similarity transformation and a translation, to the form
$(\tilde\Om_i, \III)$, where $\tilde\Om_i$ is a standard cavity.
Denote $\lam_i = l_i/|\pl(\text{conv}B)|$; one has $\sum_i \lam_i =
1$ and $\RRRR_i = -\lam_i \FFF(\tilde\Om_i)$, therefore
\begin{equation}\label{6}
\RRRR(B) = -|\pl(\text{conv}B)| \cdot (\lam_0 + \sum_{i \ne 0}
\lam_i \FFF(\tilde\Om_i)) \cdot e_2.
\end{equation}

Let the set $B_\ve$ be given by $r \le 1 - \ve f(\theta/\ve)$ in
polar coordinates $r$,\, $\theta$, where $f$ is a 1-periodic
continuous piecewise differentiable nonnegative function and $\ve$
divides $2\pi$. Suppose that $f(0) = f(1) = 0$. As $\ve \to 0$,\,
$\RRRR(B_\ve)$ tends to $-2\pi \cdot \FFF(\Om_f)$, where
\begin{equation*}\label{6+}
\Om_f = \{ (x_1, x_2):\ 0 \le x_1 \le 1,\ 0 \le x_2 \le f (x_1) \}.
\end{equation*}
Thus, the problem (\~A) (which was not yet rigorously formulated)
can be stated as follows:
 \vspace{2mm}

(A)\, {\it Find $\sup_{f} \FFF(\Om_f)$ over all continuous piecewise
differentiable nonnegative functions $f:\ [0,\, 1] \to \RRR_+$ such
that $f(0) = f(1) = 0$.}
 \vspace{2mm}

The problem (\~B) reads as:\, find $\sup_{B \subset K} |\RRRR(B)|$,
where $K \subset \RRR^2$ is a convex bounded set with nonempty
interior. In view of (\ref{6}), it amounts to the problem
 \vspace{2mm}

(B)\,  {\it Find $\sup_{\Om} \FFF(\Om)$ over all standard cavities
$\Om$.}
 \vspace{2mm}

Indeed, let $\Om_n$ be a sequence of sets solving the problem (B);
then a sequence of bodies $B_n \subset K$ approximating $K$ solves
the problem (\~B), if all the cavities of $B_n$ are similar to
$\Om_n$ and the length of the convex part of $B_n$ tends to zero.
Thus, one has $\sup_{B \subset K} |\RRRR(B)| = |\pl K| \cdot
\sup_{\Om} \FFF(\Om)$.

\section{Problem (A): numerical results}

Introduce the shorthand notation $\FFF[f] := \FFF(\Om_f)$. Note that
the functional $\FFF$ is continuous in the $C^1$ topology: if $f_n
\stackrel{C^1}{\longrightarrow} f$ then $\FFF[f_n] \to \FFF[f]$;
thus there exists a sequence of piecewise linear functions
maximizing $\FFF$. Therefore it seems natural to look for the
maximum in classes of continuous functions $f$ with piecewise
constant derivative $f'$. We also examined classes of continuous
functions with piecewise constant second derivative $f''$. In the
first case the graph of $f$ is a broken line, and in the second, a
curve composed of arcs of parabolas.

In each numerical experiment there were made $N_1 N_2$ trials with a
billiard particle in $\Om_f$. Usually $N_1$ and $N_2$ were taken
equal and varied from several hundreds to several thousands (up to
$5000$). Initially, the particle is located at $(\xi_i, 0) \in
\III$, where $\xi_i = (i - 1/2)/{N_1}$, and has the velocity
$(\sin\vphi_j,\, \cos\vphi_j)$, where $\vphi_j = \pi(j - \frac{N_2 +
1}{2})/{N_2}$,\, $i = 1,\ldots, N_1$,\, $j = 1,\ldots, N_2$. Then
the least time instant, when the particle gets into $\III$ again, is
fixed, and the particle velocity $v_{ij} = -(\sin\vphi_{ij}^+,\,
\cos\vphi_{ij}^+)$ just before this instant is registered. The sum
$$
\FFF = \frac{3}{8}\ \frac{\pi}{N_1 N_2} \sum_{i=1}^{N_1}
\sum_{j=1}^{N_2} \cos\vphi_j \left( 1 + \cos(\vphi_j - \vphi_{ij}^+)
\right)
$$
is considered to be an approximation for the integral $\FFF[f]$.

The algorithm simulating the billiard dynamics, as well as numerical
integration, were implemented in Programming Language C. The
precision accuracy achieved $10^{-6}$; it was controlled, firstly,
by the differences between the successive approximations of $\FFF$
as $N_1$ and $N_2$ increase, and secondly, by comparison with the
analytic results. To maximize resistance, there were used
optimization algorithms from the {\it Genetic Algorithm and Direct
Search Toolbox} of the computational system MATLAB, version 7.2;
these methods do not require any information of the derivatives of
the objective function.

The obtained results are as follows.
 \vspace{1mm}

1.~ In the class of two-segment broken lines
$$
f_{\al,\bt}(\xi) = \left\{
 \begin{array}{lll}
 \al\, \xi, & \text{if} & 0 \le \xi \le \xi_0\\
 \bt\, (1 - \xi), & \text{if} & \xi_0 \le \xi \le 1\,,
 \end{array}
\right.
$$
where $\al > 0$,\, $\bt > 0$,\, $0 < \xi_0 < 1$,\, $\al\, \xi_0 =
\bt\, (1 - \xi_0)$, the maximum of $\FFF$ equals $1.42621$ and is
achieved at $\al = \bt = \al_0 \approx 1.12$. Then $\xi_0 = 0.5$;
the corresponding set $\Om_{f_{\al_0,\al_0}}$,\, is an isosceles
triangle with the angle $83.6^0$ at the top vertex. It is shown on
Fig.\,(a), with the lateral sides drawn dashed.

The function $\FFF[f_{\al,\al}]$ oscillates and goes to $1.25$ as
$\al \to +\infty$ (see Fig.\,(b)).
\noindent
\begin{center}
\begin{tabular}{c c}
\includegraphics[width=0.47\columnwidth]{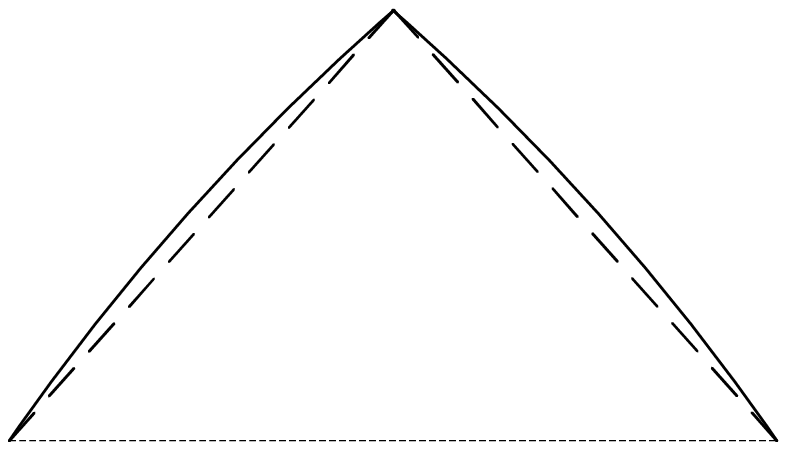}
&
\includegraphics[width=0.47\columnwidth]{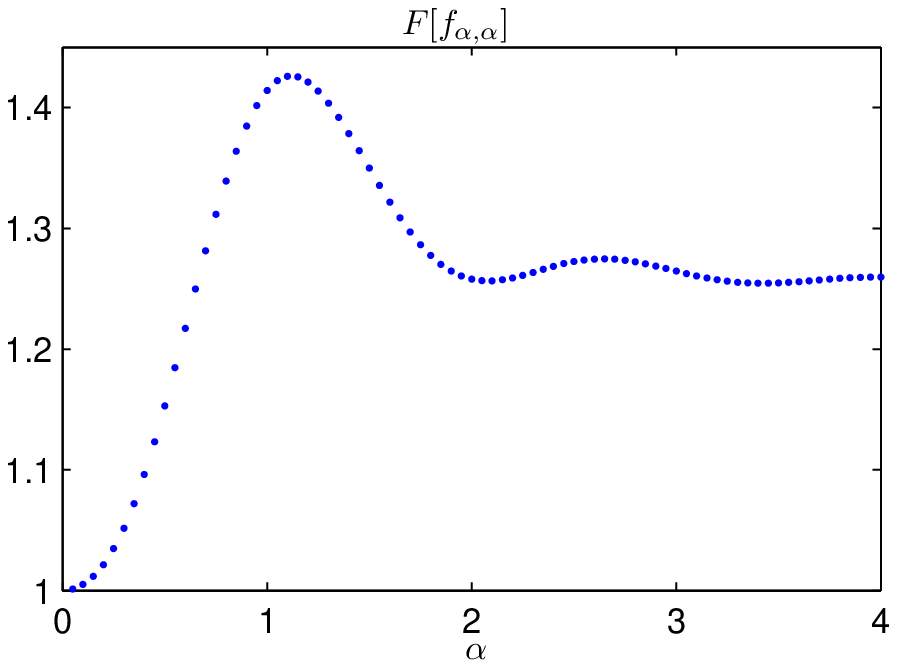} \\
(a)&(b)\\
\end{tabular}
\end{center}
 \vspace{4mm}

 2.~ In the class of two-segment piecewise quadratic functions
$$
f_{\al_1\al_2\bt_1\bt_2}(\xi) = \left\{
 \begin{array}{lll}
 \al_1 \xi^2 + \bt_1 \xi, & \text{if} & 0 \le \xi \le \xi_0\\
 \al_2 (1 - \xi)^2 + \bt_2 (1 - \xi), & \text{if} & \xi_0 \le \xi \le 1\,,
 \end{array}
\right.
$$
where $0 < \xi_0 < 1$,\, $\al_1 \xi_0^2 + \bt_1 \xi_0 = \al_2 (1 -
\xi_0)^2 + \bt_2 (1 - \xi_0)$, the maximum of $\FFF$ is achieved at
$\xi_0 = 0.5$,\, $\al_1 = \al_2 = \al \approx -0.486$,\, $\bt_1 =
\bt_2 = \bt \approx 1.361$, and is equal to $1.43816$. The
corresponding set $\Om_{f_{\al\al\bt\bt}}$ is a curvilinear
isosceles triangle; it is shown on Fig.\,(a) above with lateral
sides drawn with solid lines. Its height is equal to the height of
the optimal triangle from the item 1; so to say, this triangle is
obtained from the previous one by a slight ``bending outwards''{} its
lateral sides.
 \vspace{1mm}

3.~ In the class of broken lines with many segments the simulations
become more cumbersome. Let $x^0 = (x_1^0, x_2^0) = (0,0)$,\, $x^1 =
(x_1^1, x_2^1), \ldots, x^m = (x_1^m, x_2^m) = (1,0)$ be the
vertices of the broken line, with $m$ being the number of segments.
Experiments with relatively small $m$ $(m \le 5)$ showed that making
the broken line symmetric with respect to the vertical line $x_1 =
1/2$ and taking the values $x_1^i$,\, $i = 1, \ldots, m-1$ equally
spaced in $[0,\, 1]$ favor some moderate increase of resistance.
Therefore the posterior study was restricted to symmetric broken
lines with $x_1^i = i/m$, and thus the number of effective
parameters was reduced almost four times: from $2(m - 1)$ to
$\lfloor \frac{m}{2} \rfloor$. This decision allowed one to take a
relatively large number of segments, $1 \le m \le 18$.

There were found many ``zigzag''{} shapes providing nearly maximal
values $\FFF \approx 1.446...$; some typical shapes are shown on the
figures (a,b,c) below. The greatest found value of $\FFF$
corresponds to the shape shown on Fig.\,(b) (m = 10) and equals
$1.446227$.
\begin{center}
\begin{tabular}{c c}
\multicolumn{2}{c}{\includegraphics[width=0.4\columnwidth]{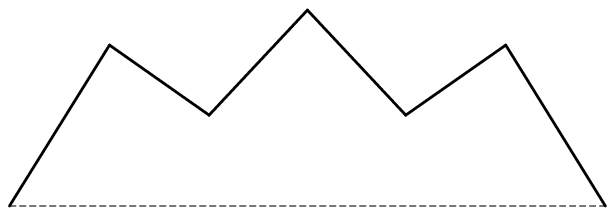}}\\
\multicolumn{2}{c}{(a) $m=6$} \vspace{0.5cm}
\\
\includegraphics[width=0.4\columnwidth]{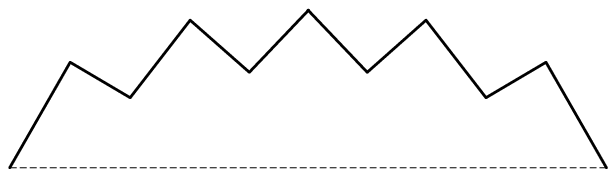}
&
\includegraphics[width=0.4\columnwidth]{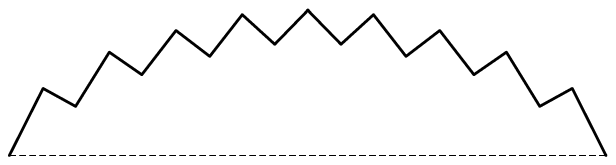} \\
(b) $m=10$&(c) $m=18$\\
\end{tabular}
\end{center}
 \vspace{4mm}

In order to verify the simulation results, as well as to find out
the common analytic form of the obtained curves, we examined the
following broken line. Fix $\Psi \in [0,\, \pi/2]$ and consider the
arc of angular size $2\Psi$ contained in the upper half-plane $x_2
\ge 0$, with the endpoints $(0,0)$ and $(1,0)$. Let $m$ be even.
Mark the points $x^0 = (x_1^0, x_2^0) = (0,0)$,\, $x^2 = (x_1^2,
x_2^2), \ldots,\, x^{2i} = (x_1^{2i}, x_2^{2i}),\, \ldots, x^m =
(x_1^m, x_2^m) = (1,0)$ on the arc, with $0 = x_1^0 < x_1^2 < \ldots
< x_1^m = 1$. Let $\del$ be the maximum of values $x_1^{2i} -
x_1^{2i-2}$. We say that a right triangle $ABC$ is {\it canonical},
if the vertex $B$ is situated above the hypotenuse $AC$ and the
median drawn from $B$ to $AC$ is vertical. For $i = 1, \ldots, m/2$,
draw the canonical triangle $\triangle x^{2i-2} x^{2i-1} x^{2i}$
with the hypotenuse $[x^{2i-2},\, x^{2i}]$. Thus, one has
$x_1^{2i-1} = \frac{1}{2}\, (x_1^{2i-2} + x_1^{2i})$,\, $x_2^{2i-1}
= \frac{1}{2}\, (x_2^{2i-2} + x_2^{2i}) + \frac{1}{2}
\sqrt{(x_1^{2i} - x_1^{2i-2})^2 + (x_2^{2i} - x_2^{2i-2})^2}$. The
broken line $x^0 x^1 \ldots x^{m-1} x^m$ composed of legs of all
triangles obtained this way will also be called {\it canonical} (the
graphs shown on the figures (a,b,c) are good approximations for
canonical lines with $m = 6,\ 10,\ 18$). If $\Psi$ is fixed and
$\del$ goes to zero, the corresponding value of $\FFF$ tends to
\begin{equation}\label{8}
\FFF(\Psi) = 1 + \frac{1}{6}\, \sin^2 \Psi + \frac{2\sqrt{2}\,
\sin\frac{\Psi}{2} - 2\, \sin^4 \frac{\Psi}{2} - \Psi}{\sin\Psi}
\end{equation}
(here $\Psi$ is expressed in radians). The proof of this convergence
is put in Appendix 2. The maximal value of $\FFF(\Psi)$ is achieved
at $\Psi_0 \approx 0.6835 \approx 39.16^0$ and is equal to
$\FFF(\Psi_0) = 1.445209$.

Let $\FFF_m(\Psi)$ be the values of $\FFF$ related to the canonical
lines with $x_1^i = i/m$,\, $i = 0, \ldots, m$. These values were
numerically calculated for various values of $\Psi$ and for $m =
6,\, 10,\, 18$. The resulting functions $\FFF_m(\Psi)$ and the
function $\FFF(\Psi)$ (solid line) are shown on the figure below.

\begin{center}
\includegraphics[width=0.65\columnwidth]{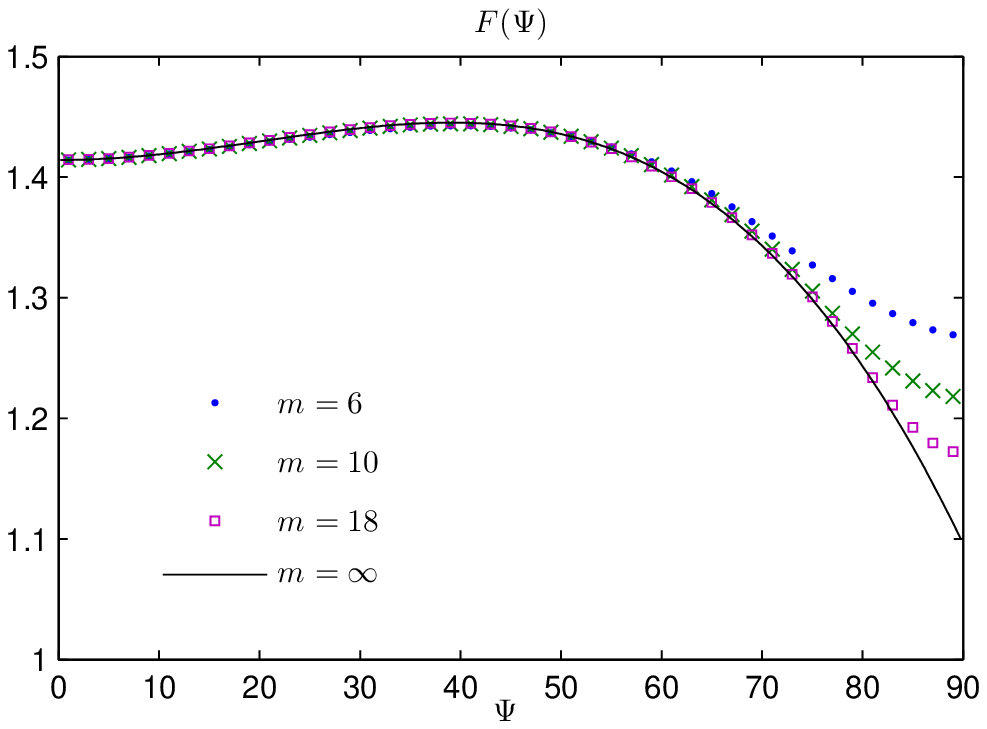}
\end{center}

Note in passing that in the limit $\Psi \to 0$ the triangles of the
corresponding canonical line approach a right isosceles triangle,
therefore the corresponding value of $\FFF$ tends to $\sqrt{2}$;
that is, $\FFF(0) = \sqrt 2$.

4.~ In the class of piecewise quadratic functions, with $\le 18$
segments, the greatest found value of $\FFF$ is 1.44772; this is the
maximal value found numerically up to the moment. The corresponding
curve is a 14-segment ``zigzag''{} line with slightly concave
segments; its visual perception is the same as of the curves shown
above on Figs.(a,b,c).
\begin{center}
\includegraphics[width=0.5\columnwidth]{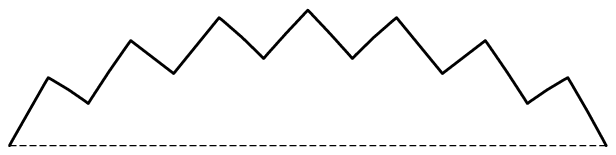}
\end{center}
Thus, one has the following estimation for
the Problem (A):
$$
1.44772 \le \sup_{f} \FFF(\Om_f) \le 1.5.
$$

5.~ We examined numerically shapes formed by infinitely small
canonical triangles. In these experiments we had to substitute the
usual billiard dynamics with the pseudo-billiard one described in
Appendix 2. Also, we tried shapes formed by various kinds of
non-canonical triangles. However, we could not increase resistance
this way.

\section{Solution of Problem (B)}

Consider the standard cavity $\Om_\ve = \Om_\ve' \cup \Om_\ve''$,
where $\Om_\ve' = [0,\, 1] \times [0,\, \ve]$ is a rectangle and
$\Om_\ve''$ is the upper semi-ellipse with the foci $F_1 = (0, \ve)$
and $F_2 = (1, \ve)$ and with the major semi-axis of length $1/\ve$.
This figure is a {\it mushroom}, with the stem $\Om_\ve'$ and the
cap $\Om_\ve''$. Note that mushroom was first proposed by Bunimovich
as an example of billiard with divided phase space \cite{Bun}.
 \vspace{45mm}

 \rput(6.5,0){
  \psarc[linewidth=0.4pt](0,0){4}{0}{180}
  \psline[linewidth=0.4pt](-4,0)(-1,0)(-1,-0.5)
  \psline[linewidth=0.4pt](4,0)(1,0)(1,-0.5)
  \psline[linewidth=0.4pt,linestyle=dashed](-1,-0.5)(1,-0.5)
  \psline[arrows=->,arrowscale=1.5](-0.7,-0.5)(-1,-0.2)(1,1.8)
      \psline[arrows=->,arrowscale=1.5](-0.7,-0.5)(-1,-0.2)
  \psline[arrows=->,arrowscale=1.5](1,1.8)(2.38,3.18)(1.2,0.7333)
  \psline(1.2,0.7333)(0.6,-0.5)
  \psdots[dotsize=2pt](-1,0)(1,0)
  \rput(-1.2,-0.17){\scriptsize $F_1$}
  \rput(1.2,-0.17){\scriptsize $F_2$}
  \rput(-4.8,1.7){\Huge $\Om_\ve$}

  }
   \vspace{17mm}

Consider the billiard particles starting from $\III = [0,\, 1]
\times \{0\}$. These particles, except for a small part of them,
make a unique reflection from the elliptical arc and then return to
$\III$, the angle between the initial and final velocity being less
that $2\arctan(\ve/2)$. The rest of the particles, i.e. those that
make reflections from the vertical sides of the stem, have the total
measure $O(\ve)$. Indeed, one easily calculates that the measure of
the particles having the first reflection from the stem is
$$
2 \cdot \int_0^1 d\xi \int_{\arctan(\xi/\ve)}^{\pi/2} \cos\vphi\,
d\vphi\ =\ 2(1 + \ve - \sqrt{1 + \ve^2}).
$$
The measure of the particles having at least one reflection from a
stem size is at most double this value.

Taking into account that the total measure of all particles incident
on $\III$ is $\int_{-\pi/2}^{\pi/2} \int_0^1 \cos\vphi\, d\xi d\vphi
= 2$ and using (\ref{5}), one gets
$$
\FFF(\Om_\ve) \ge \frac 38 \left( 2 - 4(1 + \ve - \sqrt{1 + \ve^2})
\right) \left( 1 + \cos(2\arctan\frac{\ve}{2}) \right) = 1.5 +
O(\ve).
$$
Thus, problem (B) is solved:~ $\sup_\Om \FFF(\Om) = 1.5$.
 \vspace{2mm}

Now let $K$ be a convex bounded body with nonempty interior.
Approximate it by a convex polygon $K_\ve \subset K$,\, $\ve > 0$
such that $|\pl K| - |\pl K_\ve| < \ve$. To each side $a,\ b,
\ldots$ of $K_\ve$ assign a rectangle $R_a$,\, $R_b, \ldots$ such
that one side of the rectangle (denote it by $a',\ b', \ldots$,
respectively) belongs to $a,\ b, \ldots$, etc;\, all the rectangles
belong to $K_\ve$ and do not mutually intersect;\, and the common
length of the part of perimeter of $K_\ve$ not occupied by $a',\ b',
\ldots$, is less than $\ve$.
  \vspace{60mm}

 \rput(3,0){
 \scalebox{0.8}{
 \psecurve[linewidth=0.3pt](1,0.6)(0,2)(0,3.5)(0.5,4.9)(1.5,6.1)(3,6.5)(5,6.2)(6.5,5.15)
 (7.55,3.8)(8.2,1.5)(7.8,0.75)(6,-0.05)(4,-0.2)(1,0.6)(0,2)(0,3.5)
 \psline[linewidth=0.8pt](1,0.6)(6,-0.05)
    \pspolygon[linewidth=0.9pt,linestyle=dashed](1.5,0.535)(5.75,-0.0175)(5.815,0.4825)(1.565,1.035)
    \rput(3,0.15){$a$}
    \rput(4.2,0.45){\large $R_a$}
 \psline[linewidth=0.8pt](6,-0.05)(8.2,1.5)
    \pspolygon[linewidth=0.5pt,linestyle=dashed](6.22,0.105)(7.98,1.345)(7.769,1.645)(6.009,0.405)
     \rput(7.4,0.65){$b$}
    \rput(6.7,0.9){\large $R_b$}
 \psline[linewidth=0.8pt](8.2,1.5)(6.5,5.15)
   \pspolygon[linewidth=0.5pt,linestyle=dashed](7.635,1.7594)(8,1.9294)(6.67,4.785)(6.305,4.615)
 \psline[linewidth=0.8pt](6.5,5.15)(3,6.5)
 \pspolygon[linewidth=0.5pt,linestyle=dashed](6.256,4.785)(3.196,5.965)(3.35,6.365)(6.41,5.185)
 \psline[linewidth=0.8pt](3,6.5)(0.5,4.9)
  \pspolygon[linewidth=0.5pt,linestyle=dashed](0.75,5.06)(1.006,4.66)(3.056,5.972)(2.8,6.372)
 \psline[linewidth=0.8pt](0.5,4.9)(1,0.6)
  \pspolygon[linewidth=0.5pt,linestyle=dashed](0.942,1.1)(0.558,4.4)(1.058,4.458)(1.442,1.158)
 }}
 \rput(2,2){\Huge $K$}
 \rput(6.4,2.2){\Huge $K_\ve$}

  \vspace{13mm}

On each rectangle plant out a ``seedlings of mushrooms'', as shown on
the figure (a) below for $R_a$. The sides of $R_a$ that do not
belong to $a$ are shown dashed. On Fig.\,(a), there is shown $R_a
\setminus$\,({\it union of mushrooms}), all the mushrooms being of
equal size and similar to $\Om_\ve$. The total length of the lower
horizontal part (l.h.p.) of the boundary of the obtained figure is
$1 - \ve$ times the length of the corresponding size of $R_a$, that
is, $(1 - \ve) |a'|$.

Now, plant out a ``seedlings''{} of ``mushrooms of the second order''{}
(see Fig.\,(b) below). Here the length of l.h.p. of the boundary is
$(1 - \ve)^2 |a'|$.

 \vspace{15mm}

 \rput(0,0.5){
 \scalebox{0.8}{
  \scalebox{.525}[.5]{ 
        \psline[linestyle=dashed](-0.1,-3.7)(-0.1,2)(31.85,2)(31.85,-3.7)
  \multirput(0,0)(6.4,0){5}{
 \psarc[linewidth=1.5pt](3,-3){3}{0}{180}
 \psline[linewidth=1.5pt](0,-3)(2,-3)
  \psline[linewidth=1.5pt](4,-3)(6,-3)
   \psline[linewidth=1.5pt](4,-3)(4,-3.7)
        \psline[linewidth=1.5pt](2,-3)(2,-3.7)
    \psline[linewidth=1.5pt](-0.2,-3.7)(2,-3.7)
      \psline[linewidth=1.5pt](4,-3.7)(6.2,-3.7)
  }}}}
  \rput(0.2,-1.5){(a)}
  \rput(-.3,0){\large $R_a$}

 \vspace{25mm}

  \vspace*{12mm}

 \rput(-3.5,0.5){
 \scalebox{.525}[.5]{  
        \psline[linestyle=dashed](6.3,-3.7)(6.3,2)(31.85,2)(31.85,-3.7)
  \multirput(0,0)(6.4,0){5}{
 \psarc[linewidth=1.5pt](3,-3){3}{0}{180}
 \psline[linewidth=1.5pt](0,-3)(2,-3)
  \psline[linewidth=1.5pt](4,-3)(6,-3)
    \psline[linewidth=1.5pt](4,-3)(4,-3.7)
        \psline[linewidth=1.5pt](2,-3)(2,-3.7)
         \psline[linewidth=1.5pt](4,-3.7)(8.4,-3.7)
       \rput(4.15,-3.14){
 \scalebox{.1575}[.15]{ 
  \multirput(0,0)(6.4,0){4}{
 \psarc[linewidth=10pt](3,-3){3}{0}{180}
 \psline[linewidth=10pt](0,-3)(2,-3)
  \psline[linewidth=10pt](4,-3)(6,-3)
    \psline[linewidth=10pt](4,-3)(4,-3.7)
        \psline[linewidth=10pt](2,-3)(2,-3.7)
    \psline[linewidth=10pt](-0.2,-3.7)(2,-3.7)
      \psline[linewidth=10pt](4,-3.7)(6.2,-3.7)
            \psline[linewidth=15pt,linecolor=white](2.3,-3.7)(3.7,-3.7)
  }}}
      }}}

 \rput(-2.2,-1){\pspolygon[linecolor=white,fillstyle=solid,fillcolor=white](-2,0)(-2,2.2)(2.075,2.2)(2.08,0)}
 \rput(14.85,-1){\pspolygon[linecolor=white,fillstyle=solid,fillcolor=white](-1.48,0)(-1.47,2.2)(1.5,2.2)(1.5,0)}
 \rput(0,-1.5){(b)}

 \vspace{25mm}

Continuing this process, one finally obtains the figure $\tilde R_a$
such that the length of the l.h.p. of its boundary is less than $\ve
|a'|$. Similarly, one obtains the figures $\tilde R_b, \ldots$.

Put $B_\ve = \left( K_\ve \setminus (R_a \cup R_b \cup
\ldots)\right) \cup (\tilde R_a \cup \tilde R_b \cup \ldots)$. All
the cavities of $B_\ve$ and similar to $\Om_\ve$, the length of the
convex part of $\pl B_\ve$ tends to zero and $|\pl (\text{conv}
B_\ve)| \to |\pl K|$ as $\ve \to 0$. (Recall that the convex part of
$\pl B_\ve$ is $\pl B_\ve \cap \pl (\text{conv} B_\ve)$). Therefore
$\lim_{\ve\to 0} |\RRRR(B_\ve)| = 1.5\, |\RRRR(K)|$; this solves the
maximal resistance problem (\~B).

\section{Conclusions}

There are many technical devices utilizing wind pressure force: for
example, ship sail, windmill arm, etc. We are interested here in
maximizing the pressure force of the wind consisting of {\it
non-interacting} particles. An example of such a wind is provided by
the flow of solar photons incident on a solar sail. We first studied
the problem numerically and in a restricted class of bodies, and
then, basing on the gained experience and intuition, found the
solution analytically.

On the picture below, there are shown ``pre-optimal''{} bodies, that
is, elements of maximizing sequences of bodies for problems (\~A)
and (\~B). (In the case (\~B), $K$ is the unit circle.)
 \vspace{0mm}
\begin{center}
\begin{tabular}{c c}
\includegraphics[width=0.50\columnwidth]{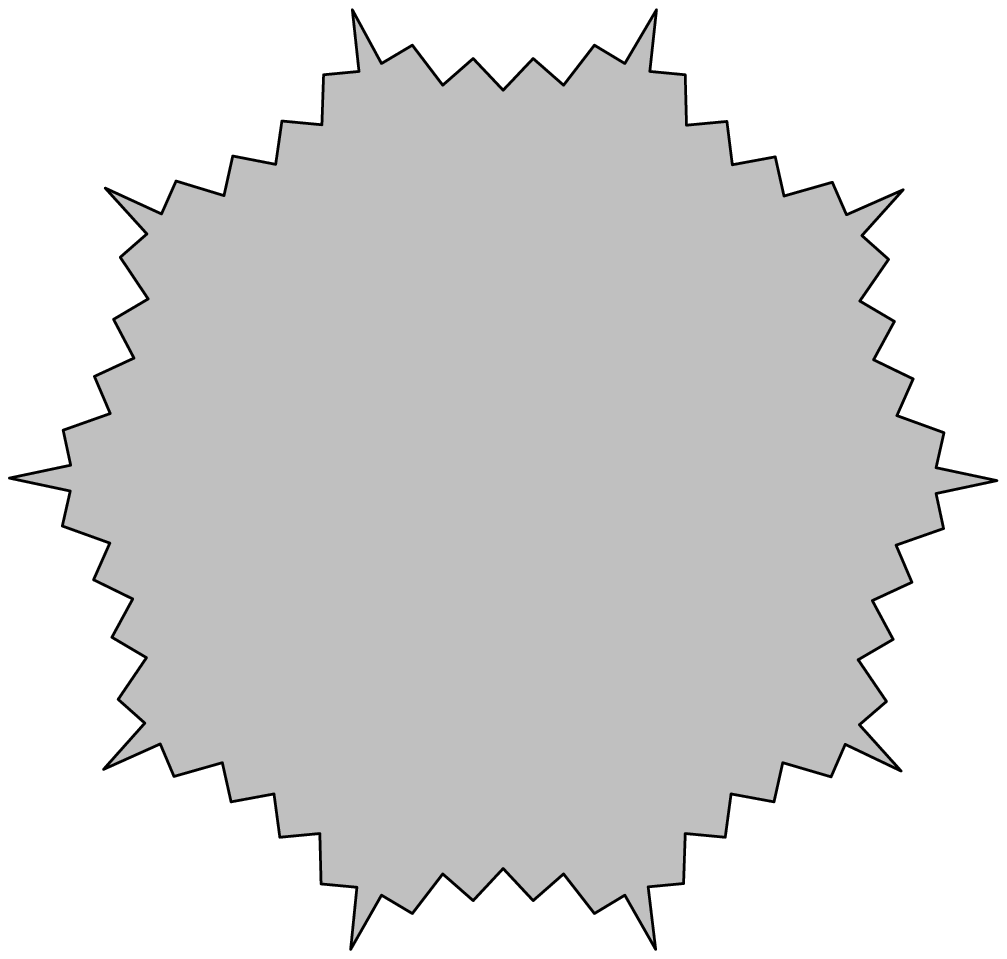}
&
\includegraphics[width=0.50\columnwidth]{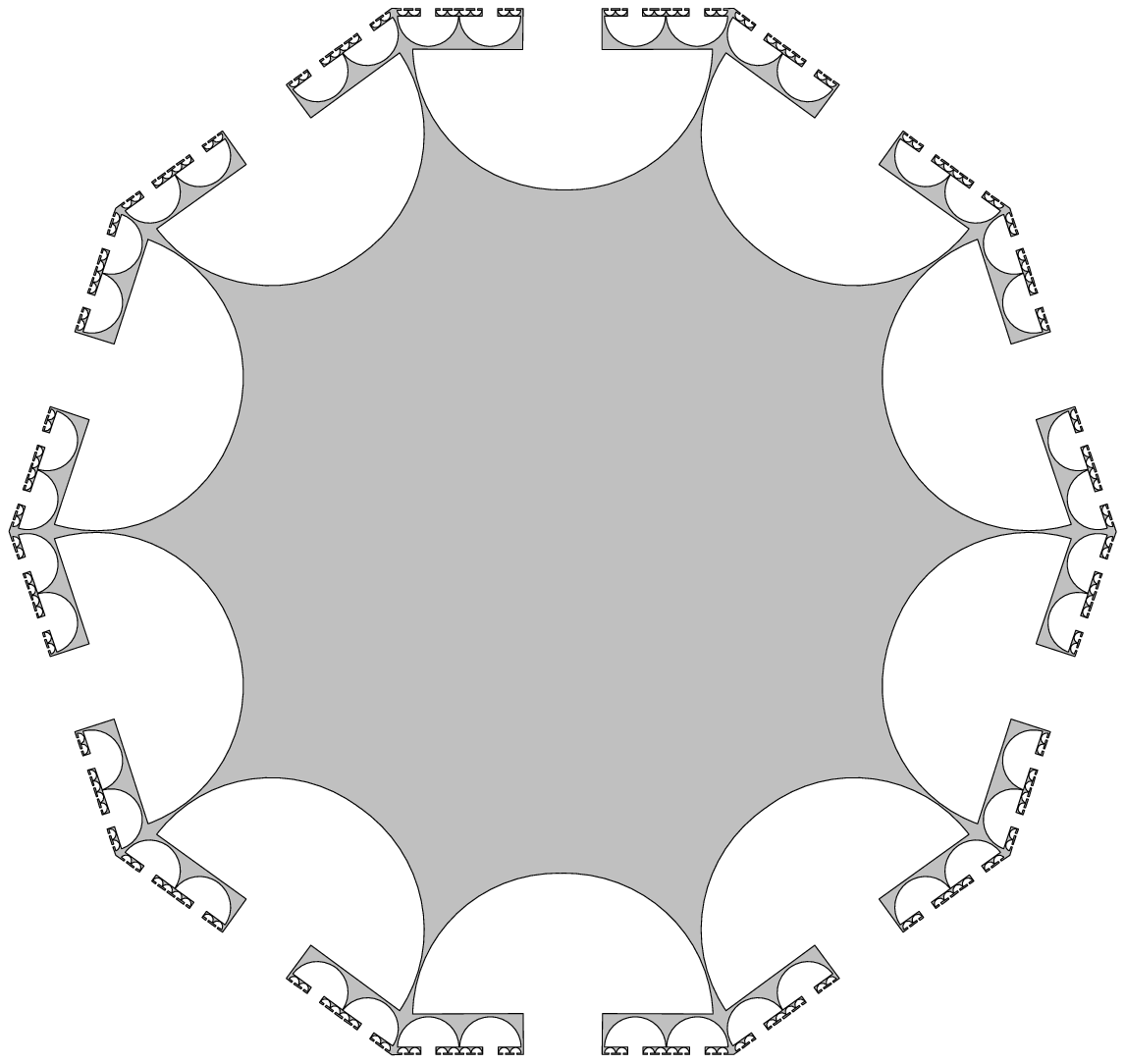} \\
(a)&(b)\\
\end{tabular}
\end{center}
 \vspace{3mm}

One can see that the boundary of the second figure is much more
complicated than of the first one. ``Complexity''{} of the boundary
can be measured by the maximal value of rotation angle of the normal
vector on a small part of the boundary. This ``maximal instantaneous
rotation''{} is, approximately, $90^0 + \Psi_0 \approx 129^0$ for
the first figure and $360^0$ for the second one.

To be precise, fix a convex bounded body $K$ and consider a sequence
$B_n$,\, $n = 1,\, 2,\ldots$ approximating $K$. Define the value
rot$(\{ B_n \})$ measuring boundary complexity of the sequence in
the following way. For each point $x \in \pl B_n$, let $\nu^n_x$ be
the outer unit normal at $x$ and let $\theta_n(x) := \arg(\nu^n_x)$.
For any $x \in \pl K$, define rot$_n(x,\ve) := \sup_{x_1,x_2 \in \pl
B_n \cap \NNN_\ve(x)} (\theta_n(x_2) -
\theta_n(x_1))$;\footnote{Recall that $\NNN_\ve(x)$ is the
$\ve$-neighborhood of $x$.} let us call it $(x,\ve)$-rotation of the
boundary $\pl B_n$. It is monotone non-increasing as $\ve \to 0$.
Then put rot$_n(\ve) = \sup_{x \in \pl K}$\,rot$_n(x,\ve)$ and
define rot$(\ve) = \lim\sup_{n\to\infty}$\,rot$_n(\ve)$:\,
$\ve$-rotation for the approximating sequence $B_n$. Finally, define
rot$(\{ B_n \}) = \lim_{\ve \to 0}\,$rot$(\ve)$ and call this value
{\it rotation of boundary} for the given sequence of bodies. This
value is $\pi/2 + \Psi_0$ for the first sequence of figures, and
$2\pi$ for the second one. The difference is almost threefold.

Further, the boundary length for the second sequence of figures
tends to infinity. On the other hand, one easily calculates that the
boundary length tends to $\frac{\sqrt{2}}{\cos(\Psi_0/2)}\, 2\pi =
1.501 \cdot 2\pi$ for the first sequence; that is, the limit value is
approximately 1.5 times the perimeter of unit circle.

Throughout this paper it was assumed that the ``inner temperature''{}
of the wind is zero. However, one can show that in the case of
positive temperature, that is, chaotic relative motion of wind
particles, the functional to be minimized is proportional to the
functional $\FFF$ (\ref{5}). Therefore all the results of this paper
remain valid in the case of positive temperature; the detailed
explication is postponed to a future paper. We are also going to
study the three-dimensional case in a future paper.

\section*{Appendix 1}

Denote by $B^{(\ve)}$ the set depicted on Fig.\,(a) in Introduction.
Here $\ve$ divides $\pi$ and conv$B^{(\ve)}$ is a $(\pi/\ve)$-sided
regular polygon inscribed in the unit circle $K_1$. Its perimeter is
$|\pl(\text{conv}B^{(\ve)})| = 2\pi\, \sin\ve/\ve$. The convex part
of $\pl B^{(\ve)}$ is a finite collection of points and all the
cavities of $B^{(\ve)}$ are isosceles right triangles similar to
$\Om = \{ (x_1,x_2): \ 0 \le x_1 \le 1, \ 0 \le x_2 \le 1/2 - |1/2 -
x_1| \}$. Applying (\ref{6}), one gets $\RRRR(B^{(\ve)}) =
-|\pl(\text{conv}B^{(\ve)})| \cdot \FFF(\Om) \cdot e_2 = -2\pi\,
(\sin\ve/\ve)\, \FFF(\Om) \cdot e_2$. According to the same formula
(\ref{6}), resistance of the unit circle $K_1$ equals $\RRRR(K_1) =
-|\pl K_1| \cdot e_2 = -2\pi \cdot e_2$. It remains to calculate
$\FFF(\Om)$.

Let $A_1$,\, $A_2$, and $A_{12}$ be the subsets of $\interval \times
[0,\, 1]$ given by the inequalities $\xi < - \tan\vphi$,\, $\xi
> 1 - \tan\vphi$, and $-\tan\vphi < \xi < 1 - \tan\vphi$,
respectively. One easily verifies the following.

(a) If $(\vphi, \xi) \in A_1$ then the corresponding billiard
particle makes only one reflection from the left leg of $\Om$ and
$\vphi_\Om(\vphi, \xi) = -\pi/2 - \vphi$.

(b) If $(\vphi, \xi) \in A_2$ then there is a single reflection from
the right leg of $\Om$ and $\vphi_\Om(\vphi, \xi) = \pi/2 - \vphi$.

(c) If $(\vphi, \xi) \in A_{12}$ then there is a double reflection,
either from the left and then from the right leg or vice versa, and
$\vphi_\Om(\vphi, \xi) = \vphi$.

Therefore, according to (\ref{5}), one has
$$
\FFF(\Om)\, =\, \frac{3}{8} \int\!\!\!\int_{A_1} \left( 1 +
\cos(\pi/2 + 2\vphi) \right)\, \cos\vphi\, d\vphi\, d\xi +
$$
$$
+ \frac{3}{8} \int\!\!\!\int_{A_2} \left( 1 + \cos(2\vphi - \pi/2)
\right)\, \cos\vphi\, d\vphi\, d\xi + \frac{3}{8}
\int\!\!\!\int_{A_{12}} 2\, \cos\vphi\, d\vphi\, d\xi\ =
$$
\begin{equation}\label{aster}
 =\ I + II + III.
\end{equation}
Direct calculation gives $I = II = \frac{3}{4} -
\frac{1}{2\sqrt{2}}$,\, $III = \frac{3}{2}\, (\sqrt{2} - 1)$; thus,
$\FFF(\Om) = \sqrt{2}$.

\section*{Appendix 2}

Recall that a right triangle is called canonical if (a) it is
situated above its hypotenuse and (b) the median dropped on the
hypotenuse is vertical. The angle $\al$ the hypotenuse forms with
the horizontal line, $\al \in (-\pi/2,\, \pi/2)$, is called
inclination of the triangle. Consider a particle that intersects the
hypotenuse, gets into the triangle, makes one or two reflections
from the legs, and then intersects the hypotenuse again and leaves
the triangle. Denote by $\vphi$ the angle the initial velocity $v$
forms with the vector $e_2 = (0, 1)$, and by $\vphi^+$, the angle
the final velocity $v^+$ forms with $-e_2$. Thus, one has $v =
(\sin\vphi, \cos\vphi)$ and $v^+ = -(\sin\vphi^+, \cos\vphi^+)$,
where $\vphi$ and $\vphi^+$ vary between $-\pi/2 + \al$ and $\pi/2 +
\al$.

Parametrize the hypotenuse by the parameter $\xi \in [0,\, 1]$; the
value $\xi = 0$ corresponds to the left endpoint of the hypotenuse,
and the value $\xi = 1$, to the right one. Like in Appendix 1,
denote by $A_1$ the set of values $(\vphi, \xi) \in [-\pi/2 + \al,\
\pi/2 + \al] \times [0,\, 1]$ corresponding to particles having a
single reflections from the left leg, by $A_2$, the set of values
corresponding to a single reflection from the right leg, and by
$A_{12}$, the set corresponding to particles having double
reflections. One easily finds that $A_1$ is given by the inequality
$\xi < -\frac{\sin\vphi}{\cos(\vphi - \al)}$,\, $A_2$, by the
inequality $\xi > 1 - \frac{\sin\vphi}{\cos(\vphi - \al)}$, and
$A_{12}$, by the double inequality $-\frac{\sin\vphi}{\cos(\vphi -
\al)} < \xi < 1 - \frac{\sin\vphi}{\cos(\vphi - \al)}$. Moreover,
for $(\vphi, \xi) \in A_1$ holds $\vphi^+ = -\pi/2 + \al - \vphi$,
for $(\vphi, \xi) \in A_2$ holds $\vphi^+ = \pi/2 + \al - \vphi$,
and for $(\vphi, \xi) \in A_{12}$ holds $\vphi^+ = \vphi$. On the
picture below, $\vphi < 0$,\, $\al > 0$ and $\xi_0 =
-\frac{\sin\vphi}{\cos(\vphi - \al)}$.
  \vspace{0mm}

 \vspace*{40mm}

  \rput(3,-0.5){
  \scalebox{0.6}{
  \psline[linewidth=1.2pt](0,0)(4,8)(8,6)
  \psline[linewidth=1.2pt,linestyle=dashed](0,0)(8,6)
   \psline[linewidth=0.9pt,linestyle=dotted](0,0)(1,0)
    \psarc[linewidth=0.4pt](0,0){0.6}{0}{36.87}
    \rput(0.85,0.3){\large $\al$}
    \psline[linewidth=0.8pt,arrows=->,arrowscale=2](2.833,-1)(2.167,1)
       \psline[linewidth=0.8pt](2.167,1)(1.5,3)
             \psline[linewidth=0.8pt,arrows=->,arrowscale=2](1.5,3)(10.5,6)
    \psline[linewidth=0.8pt,arrows=->,arrowscale=2](3.667,-1)(3,1)
       \psline[linewidth=0.8pt](3,1)(2,4)
             \psline[linewidth=0.8pt,arrows=->,arrowscale=2](2,4)(6.5,5.5)
              \psline[linewidth=0.8pt](6.5,5.5)(8,6)
         \psdot[dotsize=4pt](2.667,2)
          \rput(2.55,1.55){\large $\xi_0$}
    \psline[linewidth=0.8pt,arrows=->,arrowscale=2](4.5,-1)(3.833,1)
       \psline[linewidth=0.8pt](3.833,1)(2.5,5)
             \psline[linewidth=0.8pt,arrows=->,arrowscale=2](2.5,5)(5.5,6)
               \psline[linewidth=0.8pt,arrows=->,arrowscale=2](5.5,6)(7,6.5)(8,3.5)
   \psline[linewidth=0.9pt,linestyle=dotted](4.5,-1)(4.5,-0.2)
    \psarc[linewidth=0.4pt](4.5,-1){0.6}{90}{108.43}
              \rput(4.8,-0.5){\large $\vphi$}
  }}

 \vspace{15mm}

Consider the parallel beam of particles falling on the hypotenuse in
the direction $\vphi$. If $|\sin\vphi| < \cos(\vphi - \al)$, then
the portion of particles that make only one reflection equals
$|\sin\vphi|/\cos(\vphi - \al)$ and the direction of reflected
particles is $\pm\pi/2 + \al - \vphi$; one has to choose the sign
''+''{} if $\sin\vphi > 0$, and ``$-$''{} if $\sin\vphi < 0$. The rest
of particles make double reflections; the portion of these particles
is $1 - |\sin\vphi|/\cos(\vphi - \al)$, and the direction of
reflected particles is $\vphi$. If $|\sin\vphi| > \cos(\vphi - \al)$
then all the particles make a single reflection and the direction of
reflected particles is $\pm\pi/2 + \al - \vphi$.

Now, consider the arc of circumference of angular size $2\Psi$
contained in the half-plane $x_2 \ge 0$, with the endpoints $A =
(0,0)$ and $C = (1,0)$. Parametrize this arc with the parameter $\al
\in [-\Psi,\, \Psi]$; the value $\al = -\Psi$ corresponds to the
point $A$, and $\al = \Psi$, to the point $C$.
Divide it into a large number of small arcs and substitute each of
them with two legs of the corresponding canonical triangle. The
resulting broken line (shown on the figure below) defines a standard
cavity. Denote by $\del$ the maximum length of a small arc.

 \vspace{18mm}

 \vspace*{40mm}

 \rput(6,-2){
 \scalebox{1.3}{
 \psarc[linewidth=0.2pt,linestyle=dashed](0,0){5}{36}{144}
 \psline[linewidth=0.5pt](0,5)(0.26125,5.248)(0.5225,4.9725)(0.781,5.193)(1.0395,4.8905)
 (1.2922,5.08455)(1.545,4.7553)(1.78925,4.9232)(2.0335,4.5678)(2.267,4.7106)(2.5,4.33)(2.72,4.45)
 (2.939,4.045)(3.142,4.142)(3.346,3.716)(3.53,3.792)(3.716,3.346)(3.88,3.404)(4.045,2.939)
  \psline[linewidth=0.5pt](0,5)(-0.26125,5.248)(-0.5225,4.9725)(-0.781,5.193)(-1.0395,4.8905)
 (-1.2922,5.08455)(-1.545,4.7553)(-1.78925,4.9232)(-2.0335,4.5678)(-2.267,4.7106)(-2.5,4.33)(-2.72,4.45)
 (-2.939,4.045)(-3.142,4.142)(-3.346,3.716)(-3.53,3.792)(-3.716,3.346)(-3.88,3.404)(-4.045,2.939)
  \psline[linewidth=0.3pt,linestyle=dashed](-4.045,2.939)(4.045,2.939)
      \psline[linewidth=0.2pt,arrows=->,arrowscale=1.5](1.245,2.939)(2.058,4.068)
      \psline[linewidth=0.2pt](2.058,4.068)(2.383,4.52)
        \psline[linewidth=0.2pt,arrows=->,arrowscale=1.5](2.383,4.52)(0.383,4.21)
        \psline[linewidth=0.2pt](2.3,4.652)(2.136,4.626)
            \psline[linewidth=0.2pt,arrows=->,arrowscale=1.5](1.065,2.939)(1.878,4.068)
            \psline[linewidth=0.2pt](1.878,4.068)(2.3,4.652)
                \psline[linewidth=0.2pt,arrows=<-,arrowscale=1.5](1.2764,3.445)(2.136,4.626)
    \rput(-4.3,2.9){$A$}
  \rput(4.3,2.9){$C$}
  \rput(0,5.6){$B$}

 }}

 \vspace{-2mm}

For small $\del$, the scheme of billiard reflection can be
approximately substituted with the following description
(pseudo-billiard reflections from the arc $ABC$). A particle of some
mass moving  in a direction $\vphi \in (-\pi/2,\, \pi/2)$ is
reflected from the arc $ABC$. If $|\sin\vphi| < \cos(\vphi - \al)$,
it is split into two ``splinters'' of relative masses
$|\sin\vphi|/\cos(\vphi - \al)$ and $1 - |\sin\vphi|/\cos(\vphi -
\al)$. The first splinter is reflected in the direction $\pm\pi/2 +
\al - \vphi$, and the second, in the direction $\vphi$. If
$|\sin\vphi| \ge \cos(\vphi - \al)$, there is no splitting, and the
whole particle is reflected in the direction $\pm\pi/2 + \al -
\vphi$. The described dynamics will be called {\it pseudo-billiard}
one. A particle of unit mass starts moving at a point of $\III = AC$
in a direction $\vphi \in (-\pi/2,\, \pi/2)$, and after several
pseudo-billiard reflections, the resulting splinters return to
$\III$.

As a result of the described substitution of the billiard dynamics
with the pseudo-billiard one, one obtains the function
$$
\FFF(\Psi)\ =\ \frac 38\, \sum_i \int_{-\pi/2}^{\pi/2} \int_0^1
m_i(\vphi, \xi) \left( 1 + \cos(\vphi_i^+(\vphi,\xi) - \vphi)
\right) d\mu(\vphi, \xi),
$$
where $m_i = m_i(\vphi, \xi)$ are masses and $\vphi_i^+ =
\vphi_i^+(\vphi, \xi)$, final directions of the splinters resulting
from the particle with the initial data $(\vphi, \xi)$. (As we will
see later, splitting can really occur only once, after the first
reflection, therefore there are at most two splinters.) The
difference between $\FFF(\Psi)$ and the true value of the functional
$\FFF$ (\ref{5}) is $O(\del)$,\, $\del \to 0^+$. This fact can be
expressed as $\lim_{\del \to 0^+} \FFF(\text{broken line}) =
\FFF(\Psi)$. Below we will calculate $\FFF(\Psi)$.

In order to describe the pseudo-billiard motion, it is helpful to
change the variables. Consider the circumference containing the arc
under consideration, and parametrize it with the same angular
variable $\al$; this time $\al$ varies in $[-\pi,\, \pi]$. Consider
a particle that starts moving at some point $\bt$ of the
circumference, intersects $\III$ at some point $(\xi, 0)$, and then
reflects from the arc, according to the pseudo-billiard rule, at a
point $\al$. Thus, one has $\Psi < |\bt| \le \pi$ and $|\al| <
\Psi$. If $|\bt| \le \pi/2$, there is no splitting, and if $|\bt| >
\pi/2$, there is.

Let us describe the dynamics of the first splinter. For a while,
change the notation; let $\bt =: \al_{-1}$,\, $\al =: \al_0$, and
designate by $\al_1$ the point of intersection of the splinter
trajectory with the circumference. Denote by $\vphi$ the initial
direction of the particle, and by $\vphi'$, the direction of the
splinter after the first reflection. (We do not call it $\vphi^+$,
since there may be more reflections.) One has $\vphi = \pm \pi/2 +
(\al_0 + \al_{-1})/2$ and $\vphi' = \pi/2 + (\al_0 + \al_1)/2$.
Then, taking into account that $\vphi' = \pm\pi/2 + \al_0 - \vphi$,
one gets
$$
\vphi = (\al_0 - \al_1)/2 ~~~ \text{and} ~~~ (\al_{-1} + \al_1)/2 =
\pi/2,
$$
the equalities being true mod\,$\pi$. In other words, the points
$\al_{-1}$ and $\al_1$ lie on the same vertical line; see the
figure.
 \vspace{0mm}
 \newpage
 \vspace*{5mm}

 \rput(5,-2){
 \scalebox{0.8}{
 \pscircle[linewidth=0.2pt](0,0){4}
 \psline[linewidth=0.35pt,arrows=->,arrowscale=2.5](1.7,-3.62)(-3,2.646)(1.7,3.62)(-3,-2.646)
 \psline[linestyle=dashed](-3.625,1.69)(3.625,1.69)
 \psarc[linewidth=1.3pt](0,0){4}{25}{155}
 \psline[linewidth=0.5pt,linestyle=dotted](1.7,-3.62)(1.7,3.62)
 \psline[linewidth=0.5pt,linestyle=dotted](-3,2.646)(-3,-2.646)
 \psdot(-2.29,1.69)
 \rput(-2.02,1.98){\scriptsize $(\xi,0)$}
 \rput(2.35,1.4){\large $\III$}
 \rput(-4.09,1.7){\large $-\Psi$}
  \rput(3.95,1.7){\large $\Psi$}
   \rput(2.3,-3.8){$\al_{-1} = \bt$}
   \rput(-3.6,2.8){$\al_0 = \al$}
   \rput(1.9,3.8){$\al_1$}
   \rput(-3.25,-2.8){$\al_2$}
   \rput(0.3,-1.4){$\vphi$}
   \rput(-0.5,2.9){$\vphi'$}
   \rput(-2.1,-2){$\vphi^+$}

   \psline[linewidth=0.35pt,arrows=->,arrowscale=2.5](1.7,-3.62)(-0.18,-1.114)
   \psline[linewidth=0.35pt,arrows=->,arrowscale=2.5](-3,2.646)(0.1,3.289)

  }}

 \vspace{58mm}

If $\al_1$ belongs to the arc $[-\Psi,\, \Psi]$ then there occurs
one more reflection, this time without splitting, since the splinter
arrived from the point $\al_0 \in \interval$. Extend the trajectory
after the second reflection until the intersection with the
circumference at a point $\al_2$. Using an argument analogous to the
one stated above, one derives the formula $\al_0 + \al_2 = \pi$; it
follows that the point $\al_2$ does not belong to the arc, that is,
there are no reflections anymore.

Summarizing, the pseudo-billiard dynamics is as follows. After the
first reflection from the arc, the particle may, and may not, split
into two ``splinters''. If $\al_{-1} \in [-\pi/2,\ - \Psi] \cup
[\Psi,\ \pi/2]$, there are no splitting, and the reflection is
unique. If $\al_{-1} \in [-\pi,\ -\pi/2) \cup (\pi/2,\ \pi]$, there
is splitting into two splinters. If $\al_{-1} \in [-\pi + \Psi,\
-\pi/2) \cup (\pi/2,\ \pi - \Psi]$, the first splinter makes no
reflections anymore. If $\al_{-1} \in [-\pi,\ -\pi + \Psi] \cup [\pi
- \Psi,\ \pi]$, it makes one more reflection (without splitting)
from the arc, and the final direction is $\vphi^+ = \pi/2 + (\al_1 +
\al_2)/2$. Taking into account the above equalities, one gets $\vphi
- \vphi^+ = \al_{-1} + \al_0 + \pi$.

Note that the factor $1 + \cos(\vphi - \vphi^+)$, meaning the impact
force per unit mass, equals 2 for the second splinter. For the first
splinter that makes no reflections, as well as for the reflection
without splitting, this factor equals $1 + \cos(\vphi - \vphi') = 1
+ |\sin\al_{-1}|$. Finally, for the first splinter that makes one
more reflection, this factor equals $1 + \cos(\vphi - \vphi^+) = 1 -
\cos(\al_0 + \al_{-1})$.

Let us pass from the variables $\vphi$ and $\xi$ to $\al = \al_0$
and $\bt = \al_{-1}$ and calculate the integral $\FFF(\Psi)$ in
terms of the new variables. The points $\al$ and $\bt$ on the
circumference have the cartesian coordinates $\frac{1}{2\sin\Psi}\,
(\sin\Psi + \sin\al, -\cos\Psi + \cos\al)$ and
$\frac{1}{2\sin\Psi}\, (\sin\Psi + \sin\bt, -\cos\Psi + \cos\bt)$,
respectively. The interval with the endpoints $\al$ and $\bt$
intersects with the interval $\III$ at the point $(\xi, 0)$, where
\begin{equation}\label{a1}
\xi\, =\, \frac{\sin(\Psi + \al) - \sin(\Psi + \bt) + \sin(\bt -
\al)}{2\sin\Psi\, (\cos\al - \cos\bt)}\,.
\end{equation}
Further, one has
\begin{equation}\label{a2}
\vphi\, =\, \frac{\al + \bt \pm\pi}{2}\,;
\end{equation}
one has to take the sign ``$-$'' or ``$+$'', if $\bt > 0$ or $\bt < 0$,
respectively. Therefore, $\cos\vphi = |\sin \frac{\al + \bt}{2}|$.

The point $(\al, \bt)$ runs the set $[-\Psi,\, \Psi] \times \left(
[-\pi,\ -\Psi] \cup [\Psi,\ \pi] \right)$, and the mapping $(\al,
\bt) \mapsto (\vphi, \xi)$ given by (\ref{a1}),(\ref{a2}) is a
one-to-one mapping from this set to $\interval \times [0,\, 1]$,
with the Jacobian
$$
\frac{D(\vphi, \xi)}{D(\al, \bt)}\ =\ \frac{\cos(\Psi + \al) +
\cos(\Psi + \bt) - 2\cos(\bt - \al)}{4\sin\Psi\, (\cos\bt -
\cos\al)}\ + \hspace*{12mm}
$$
$$
+\ \frac{\sin(\Psi + \bt) - \sin(\Psi + \al) + \sin(\al -
\bt)}{4\sin\Psi\, (\cos\bt - \cos\al)^2}\ (\sin\al + \sin\bt)\ =
\hspace*{7mm}
$$
\begin{equation}\label{a3}
\hspace*{56mm}=\ \frac{1}{4\sin\Psi}\ \frac{\sin \frac{\al -
\bt}{2}}{\sin \frac{\al + \bt}{2}}\,;
\end{equation}
this implies that the integration factor equals
$$
\cos\vphi\, d\vphi\, d\xi\, =\, \frac{1}{4\sin\Psi}\ \Big|\sin
\frac{\bt - \al}{2} \Big|\, d\al\, d\bt.
$$
Further, the mass of the first splinter is $|\cos \frac{\al +
\bt}{2}| \Big/ |\sin \frac{\bt - \al}{2}|$, and of the second one,
$1 - |\cos \frac{\al + \bt}{2}| \Big/ |\sin \frac{\bt - \al}{2}|$.
Note also that integrating over $\bt \in [-\pi,\ -\Psi] \cup [\Psi,\
\pi]$ can be substituted with integrating over $\bt \in [\Psi,\
\pi]$ with subsequent duplication of the result. With this
substitution, one always has $\sin \frac{\bt - \al}{2} > 0$.

The integral $\FFF(\Psi)$ can be written down as the sum $\FFF(\Psi)
= I + II + III + IV$, where
$$
I\ =\ \frac{3}{16\sin\Psi}\, \int_{-\Psi}^{\Psi} d\al
\int_\Psi^{\pi/2} (1 + \sin\bt)\, \sin \frac{\bt - \al}{2}\, d\bt,
$$
$$
II\ =\ \frac{3}{16\sin\Psi}\, \int_{-\Psi}^{\Psi} d\al
\int_{\pi/2}^{\pi} 2 \left( \sin \frac{\bt - \al}{2} - \Big|\cos
\frac{\al + \bt}{2} \Big| \right) d\bt,
$$
$$
III\ =\ \frac{3}{16\sin\Psi}\, \int_{-\Psi}^{\Psi} d\al
\int_{\pi/2}^{\pi - \Psi} (1 + \sin\bt)\, \cos \frac{\al + \bt}{2}\,
d\bt,
$$
$$
IV\ =\ \frac{3}{16\sin\Psi}\, \int_{-\Psi}^{\Psi} d\al \int_{\pi -
\Psi}^{\pi} (1 - \cos(\al + \bt))\, \Big| \cos \frac{\al + \bt}{2}
\Big|\, d\bt.
$$
As a result of simple calculation, one obtains
$$
I = III = \frac{3}{16\sin\Psi}\, \left[ 4\sin\Psi  -
\frac{8\sqrt{2}}{3}\, \sin\frac{\Psi}{2} - \frac{16}{3}\, \sin^4
\frac{\Psi}{2} \right],
$$
$$
II = \frac{3}{16\sin\Psi}\, \left[ 16\, \sqrt{2}\,
\sin\frac{\Psi}{2} - 8\, \Psi \right],
$$
$$
IV = \frac{3}{16\sin\Psi}\, \left[ - \frac{8}{3}\, \sin\Psi +
\frac{8}{9}\, \sin^3 \Psi + \frac{8}{3}\, \Psi \right].
$$
Summing these expressions, one finally comes to the formula
(\ref{8}):
$$
\FFF(\Psi) = 1 + \frac{1}{6}\, \sin^2 \Psi + \frac{2\sqrt{2}\,
\sin\frac{\Psi}{2} - 2\, \sin^4 \frac{\Psi}{2} - \Psi}{\sin\Psi}\,.
$$

\section*{Acknowledgements}

This work was supported by {\it Centre for Research on Optimization
and Control} (CEOC) from the ''{\it Funda\c{c}\~{a}o para a
Ci\^{e}ncia e a Tecnologia}'' (FCT), cofinanced by the European
Community Fund FEDER/POCTI.


\begin{thebibliography}{99}

\bibitem{N}
I. Newton,\, {\it Philosophiae naturalis principia mathematica}\,
1687.

\bibitem{BK} G. Buttazzo and B. Kawohl.
\textit{On Newton's problem of minimal resistance}. Math. Intell.
{\bf 15}, 7-12 (1993).

\bibitem{BrFK} F. Brock, V. Ferone, and B. Kawohl.
\textit{A symmetry problem in the calculus of variations}. Calc.
Var. {\bf 4}, 593-599 (1996).

\bibitem{BFK} G. Buttazzo, V. Ferone, and B. Kawohl.
\textit{Minimum problems over sets of concave functions and related
questions}. Math. Nachr. {\bf 173}, 71-89 (1995).

\bibitem{BG}
G. Buttazzo and P. Guasoni,\, {\it Shape optimization problems over
classes of convex domains},\, J. Convex Anal. {\bf 4}, 343-351
(1997).

\bibitem{LP1} T. Lachand-Robert and M.~A. Peletier.
\textit{Newton's problem of the body of minimal resistance in the
class of convex developable functions}. Math. Nachr. {\bf 226},
153-176 (2001).

\bibitem{LP2} T. Lachand-Robert, M.~A. Peletier.
\textit{An example of non-convex minimization and  an application to
Newton's problem of the body of least resistance}. Ann. Inst. H.
Poincar\'e, Anal. Non Lin. {\bf 18}, 179-198 (2001).

\bibitem{CL1} M. Comte and T. Lachand-Robert.
\textit{Newton's problem of the body of minimal resistance under a
single-impact assumption}. Calc. Var. Partial Differ. Equ. {\bf 12},
173-211 (2001).

\bibitem{CL2} M. Comte and T. Lachand-Robert.
\textit{Existence of minimizers for Newton's problem of the body of
minimal resistance under a single-impact assumption}. J. Anal. Math.
{\bf 83}, 313-335 (2001).

\bibitem{LO} T. Lachand-Robert and E. Oudet.
\textit{Minimizing within convex bodies using a convex hull method}.
SIAM J. Optim. {\bf 16}, 368-379 (2006).

\bibitem{P1}
A.\,Yu. Plakhov. {\it Newton's problem of a body of minimal
aerodynamic resistance}, Dokl. Akad. Nauk {\bf 390}, 314-317 (2003).

\bibitem{P2}
A.\,Yu. Plakhov. {\it Newton's problem of the body of minimal
resistance with a bounded number of collisions}, Russ. Math. Surv.
{\bf 58}, 191-192 (2003).

\bibitem{average-04}
A.\,Yu. Plakhov. {\it Newton's problem of the body of minimum mean
resistance}, Sbornik: Mathematics {\bf 195}, N$^\text{o}$7-8,
1017-1037 (2004).

\bibitem{Bun}
L.\,B. Bunimovich. {\it Mushrooms and other billiards with divided
phase space}, Chaos {\bf 11}, 802-808 (2001).

\end{thebibliography}
\end{document}